\def\[{[\![}
\def\]{]\!]}
\def\l{\langle}
\def\r{\rangle}
\def\q{\quad}
\def\e{\eqno}
\def\n{\noindent}

\def\s{\smallskip}
\def\a{\alpha}

\noindent
\centerline{\bf  A DESCRIPTION OF THE QUANTUM SUPERALGEBRA U$_q$[OSP(2N+1/2M)]}
\centerline{\bf  VIA GREEN GENERATORS }

\vskip 32pt
\noindent
\hskip 66mm  {\bf Tchavdar D. Palev}\footnote*{Permanent Address: 
 Institute for Nuclear Research and Nuclear Energy, 1784 Sofia, 
Bulgaria; $e$-mail: tpalev@inrne.acad.bg}
\bigskip
\noindent
\centerline{International School for Advanced Studies, via Beirut 2-4,
34013 Trieste, Italy} 
\smallskip
\centerline{  and  }
\smallskip
\centerline{ International Centre for Theoretical Physics, P.O. Box 586,
34100 Trieste, Italy}

\vskip 48pt
\noindent
Talk presented at the XXI International Colloquium on Group
Theoretical Methods in Physics \hfill\break
(15-20 July 1996, Goslar, Germany)

\vskip 48pt
\noindent

{\bf 1. Introduction. Outline of the results}

\bigskip
In the present talk I'll describe the orthosymplectic Lie
superalgebra $osp(2n+1/2m)$ and also its $q-$deformed analogue
$U_q[osp(2n+1/2m)]$ in terms of a new set of generators, called
Green generators. These generators are very different form the well
known Chevalley generators. Let me underline from the very beginning 
that I am not going to consider new deformation of $U_q[osp(2n+1/2m)]$.
The deformation will be the known Hopf algebra deformation as given, 
for instance, in [1-4]. The description, however, will be given in
terms of new free generators.

For me personally the interest in the construction stems from the
observation that the Green generators are of a direct physical
significance. In a certain representation of $osp(2n+1/2m)$ 
part of these generators are Bose operators, whereas the rest are
Fermi operators. Considered as elements from the universal
enveloping algebra, the Green generators are para-Bose and
para-Fermi operators [5].
To begin with I'll state the final result. It is contained in the 
following 

\s
{\bf Theorem.} {\it $U_q[osp(2n+1/2m)]$ is an associative 
superalgebra with 1, generators 
$a_i^\pm,\; L_i,\;{\bar L_i}\equiv L_i^{-1}$, $i=1,2,\ldots,m+n=N$,
relations ($\xi,\eta =\pm\; or\; \pm 1,\;\;{\bar q}\equiv q^{-1}$)
$$
\eqalign{
&  L_iL_i^{-1}=L_i^{-1}L_i=1, \q  L_iL_j=L_jL_i,\cr
& L_ia_j^\pm=q^{\pm \delta_{ij}(-1)^{\l i \r} }a_j^\pm,\cr
& \[a_i^-,a_i^+\]=-2{L_i-{\bar L}_i\over q-{\bar q}}, \cr 
& \[\[a_i^{\eta},a_{i+ \xi}^{-\eta}\],  
a_j^{\eta}\]_{q^{-\xi (-1)^{\l i \r}\delta_{ij}}}
=2(\eta)^{\l j \r}\delta _{j,i + \xi}L_j^{-\xi \eta}a_i^{\eta}, \cr  
& [[a_{N-1}^\xi,a_N^\xi],a_{N}^\xi]_{\bar q}=0 . \cr
}\eqno(2)
$$
and ${\bf Z}_2$-grading induced from
$$
{\rm deg}(L_i)={\bar 0},\quad
{\rm deg}(a_i^\pm)= \langle i \rangle \equiv \cases 
{\bar 1, & for $i\le m$ \cr
\bar 0, &  for $i> m$. \cr} \eqno(3)
$$
}
Here and throughout
$$
[x,y]_q=xy-qyx,\quad \{x,y\}_q=xy+qyx,\q 
\[x,y\]_q=xy-(-1)^{{\rm deg}(x){\rm deg}(y)}qyx.  
\eqno(4)
$$

This theorem extends the results of several previous publications.
The first deformation of one pair of para-Bose operators was given
independently in [6] and [7]. The second paper includes also all
Hopf algebra operations. This result was generalized to any number
of parabosons in [8, 9], including some representations in the
root of unity case [10]. A similar problem for any number of
parafermions was solved in [11]. The deformation of one pair of
parafermions and one pair of parabosons was carried out in [12].
Finally, the nondeformed version of the present investigation is 
given in [13].

 The plan of the exposition will be the following. First in Sect 2 I'll 
recall the definition of the orthosimplectic Lie superalgebra (LS)
$osp(2n+1/2m)$ in a matrix form. As  next steps, a description of its 
universal enveloping algebra in terms of operators, called preoscillator 
generators (Sect. 3), and via Green generators (Sect. 4) will be given. 
Finally, in Sect. 5 the deformed algbra will be considered
and some indications of how the proof of the Theorem goes will be
mentioned.

\vskip 32pt

\noindent
{\bf 2. Definition of $osp(2n+1/2m)$ in a matrix form} [14]

\bigskip
The Lie superalgebra  $osp(2n+1/2m)$ can be defined as the set of all 
$(2n+2m+1)\times (2n+2m+1)$ matrices of the form (T=transposition)

$$
  \pmatrix{a     &  b     &  u     &   x &  x_1 \cr
           c     & -a^T   &  v     &   y &  y_1 \cr
	  -v^T   & -u^T   &  0     &   z &  z_1 \cr
	   y_1^T &  x_1^T &  z_1^T &   d &  e   \cr
	  -y^T   & -x^T   & -z^T   &   f & -d^T \cr},\eqno(5)   
$$
where 
$a$ is any $n\times n$ matrix, 
$b$ and $c$ are skew symmetric $n\times n$ matrices, 
$d$ is any $m\times m$ matrix,
$e$ and $f$ are symmetric $m\times m$ matrices, 
$x,x_1,y,y_1$ are $n\times m$ matrices, 
$u$ and $v$ are $n\times 1 $ columns, 
$z,z_1$ are $1\times m$ rows.
The even subalgebra consists of all matrices with
$
x=x_1=y=y_1=z=z_1=0,
$  
namely
$$
  \pmatrix{a     &  b       &  u     &   0 &  0 \cr
           c     & -a^T   &  v     &   0 &  0   \cr
	  -v^T   & -u^T   &  0     &   0 &  0   \cr
	   0     &  0     &  0     &   d &  e   \cr
	   0     &  0     &  0     &   f & -d^T \cr},\eqno(6)   
$$
and it is isomorphic to the Lie algebra $so(2n+1)\oplus sp(2m)$. The odd
subspace is given with all matrices
$$
  \pmatrix{0     &  0     &  0     &   x &  x_1 \cr
           0     &  0     &  0     &   y &  y_1 \cr
	   0     &  0     &  0     &   z &  z_1 \cr
	   y_1^T &  x_1^T &  z_1^T &   0 &  0   \cr
	  -y^T   & -x^T   & -z^T   &   0 &  0   \cr}.\eqno(7)   
$$
The product (= the supercommutator) is defined on any two homogeneous
elements $a$ and $b$ as
$$
\[a,b\]=ab-(-1)^{deg(a)deg(b)}ba. \eqno(8)
$$

Let $L(n/m)$ be the $2(n+m)$-dimensional ${\bf Z}_2$-graded subspace,
consisting of all matrices
$$
  \pmatrix{0     &  0     &  u     &   0 &  0   \cr
           0     &  0     &  v     &   0 &  0   \cr
	  -v^T   & -u^T   &  0     &   z &  z_1 \cr
	   0     &  0     &  z_1^T &   0 &  0   \cr
	   0     &  0     & -z^T   &   0 &  0   \cr}.\eqno(9)   
$$

\n Label the rows and the columns with the indices
$
A, B=-2n,-2n+1,\ldots,-2,-1,0,1,2,\ldots,2m
$
and let $e_{AB}$ be a matrix with 1 at the intersection of the
$A^{th}$ row and the $B^{th}$ column and zero elsewhere. Then the 
following elements (matrices) constitute a basis in $L(n/m)$:
$$
\eqalign{
& a_i^-\equiv B_i^-=\sqrt{2}(e_{0,i}-e_{i+m,0}),\quad 
  a_i^+\equiv B_i^+=\sqrt{2}(e_{0,i+m}+e_{i,0}),\quad
  i=1,\ldots,m,\cr
& a_{j+m}^-\equiv F_j^-=\sqrt{2}(e_{-j,0}-e_{0,-j-n}),\quad 
  a_{j+m}^+\equiv F_j^+=\sqrt{2}(e_{0,-j}-e_{-j-n,0}),\quad
  j=1,\ldots,n,\cr
}\eqno(10)
$$
with
$
{\rm deg}(a_i^\pm)= \langle i \rangle .
$

\s
{\bf Proposition 1.} 
{\it The LS $osp(2n+1/2m)$ is generated from
$a_i^\pm,\;i=1,\ldots,m+n\equiv N.$ }

\s
It is straightforward to show that 
$$
osp(2n+1/2m)
=lin.env. \{ a_i^\xi,\; \[a_j^\eta, a_k^\varepsilon\] | 
i,j,k=1,\ldots,N,\q \xi,\eta,\varepsilon =\pm \}. \eqno(11)
$$
Hence any further supercommutator between
$ a_i^\xi,\; \[a_j^\eta, a_k^\varepsilon \],\;\;
\xi,\eta,\varepsilon =\pm, $
is a linear combination of the same type elements. A more precise
computation gives:
$$
\[ \[a_i^\xi,a_j^\eta\],a_k^\epsilon\]=
2\epsilon^{\langle k \rangle}\delta_{jk}\delta_{\epsilon,-\eta}a_i^\xi
-2\epsilon^{\langle k \rangle}
(-1)^{\langle j \rangle \langle k \rangle}\delta_{ik}
\delta_{\epsilon,-\xi}a_j^\eta . \eqno(12)
$$
Eqs. (12) are among the supercommutation relations of all
Cartan-Weyl generators
$$
 a_i^\xi,\; \[a_j^\eta, a_k^\varepsilon\], \q  
i,j,k=1,\ldots,N,\q \xi,\eta,\varepsilon =\pm .\e(13)
$$
The rest of the supercommutation relations follow from
(12) and the (graded) Jacoby identity
($ i,j,k,l=1,\ldots,N,\q \xi,\eta,\varepsilon,\varphi =\pm$):
$$
\eqalign{
& \[ \[a_i^\xi,a_j^\eta\], \[a_k^\epsilon,a_l^\varphi \]  \] 
= 2\epsilon^{\l k \r}\delta_{jk}\delta_{\epsilon,-\eta}
\[a_i^\xi,a_l^\varphi\]
-2\epsilon^{\l k \r} (-1)^{\l j \r \l k \r}\delta_{ik}
\delta_{\epsilon,-\xi} \[a_j^\eta ,a_l^\varphi \] \cr
&\cr
&  \q\; -2\varphi^{\l l \r}(-1)^{\l j \r \l k \r} 
\delta_{jl}\delta_{\varphi,-\eta}\[a_i^\xi,a_k^\epsilon\]
 +2\varphi^{\l l \r}(-1)^{\l i \r \l j \r + \l i \r \l k \r}
\delta_{il}\delta_{\varphi,-\xi}\[a_j^\eta,a_k^\epsilon\]. \cr
}\e(14)
$$

\vskip 6mm

\n{\bf 3. Description of $U[osp(2n+1/2m)]$ via preoscillator generators}

\bigskip
The relations (12) are representation independent. More precisely,
the universal enveloping algebra (UEA) $U[osp(2n+1/2m)]$ of
$osp(2n+1/2m)$ is by definition the (free) associative algebra
with 1 of the indeterminates
$
a_1^\pm,a_2^\pm,\ldots ,a_{m+n}^\pm \equiv a_N^\pm,
$
subject to the relations (12) and (14). Since however Eqs. (14)
follow from (12), we have

\bigskip
{\bf Proposition 2.} 
{\it
(1) $U[osp(2n+1/2m)]$ is the associative unital
    algebra with generators  
    $$
    a_1^\pm,a_2^\pm,\ldots,a_{m-1}^\pm,a_{m}^\pm,a_{m+1}^\pm,
    \ldots,a_{m+n}^\pm \equiv a_N^\pm,\e(15)
    $$
    relations
    $$
    \[ \[a_i^\xi,a_j^\eta\],a_k^\epsilon\]=
    2\epsilon^{\langle k \rangle}\delta_{jk}\delta_{\epsilon,-\eta}a_i^\xi
    -2\epsilon^{\langle k \rangle}
    (-1)^{\langle j \rangle \langle k \rangle}\delta_{ik}
    \delta_{\epsilon,-\xi}a_j^\eta  \eqno(16)
    $$
    and ${\bf Z}_2$-grading induced from
    $$
    {\rm deg}(a_i^\pm)= \langle i \rangle. \e(17)     
    $$
\n (2)
$$
osp(2n+1/2m)=lin.env.\{ a_i^\xi,\; \[a_j^\eta, a_k^\varepsilon\] | 
i,j,k=1,\ldots,N,\q \xi,\eta,\varepsilon =\pm \}, 
\eqno(18)
$$
with a natural supercommutator (turning every associative superalgebra 
into a Lie superalgebra):
$$
\[a,b\]=ab-(-1)^{deg(a)deg(b)}ba. \eqno(19)
$$
}

The above proposition gives a definition of 
$U[osp(2n+1/2m)]$ in terms of a new set of generators, 
which are very different from the Chevalley
generators. The relevance of the generators $a_i^\pm$ stems from the
following observation.
The operators $a_i^\pm\equiv B_i^\pm$ 
with $i=1,\ldots,m$ satisfy the triple relations
$$
[\{B_i^\xi,B_j^\eta\},B_k^\varepsilon]=
     2\varepsilon\delta_{jk}\delta_{\varepsilon,-\eta}B_i^\xi
     +2\varepsilon\delta_{ik}\delta_{\varepsilon,-\xi}B_j^\eta, \eqno(20)
$$
whereas  $a_{i+m}^\pm\equiv F_i^\pm$ with $i=1,\ldots,n$ yields: 
$$
[[F_i^\xi,F_j^\eta],F_k^\varepsilon]=
           2\delta_{jk}\delta_{\varepsilon,-\eta}F_i^\xi
          -2\delta_{ik}\delta_{\varepsilon,-\xi}F_j^\eta. \eqno(21)
$$

\smallskip
The relations (20) and (21) are known in quantum fieid theory. They are
defining relations for para-Bose and for para-Fermi creation and
annihilation operators, respectively [5].
The para-Fermi operators generate the Lie algebra
$so(2n+1)$ [15], whereas $m$ pairs of para-Bose operators
generate a Lie superalgebra [16], which is isomorphic
to $osp(1/2m)$ [17].

In the Fock representation the para-Bose (resp. the para-Fermi)
operators become usual Bose (resp. Fermi) operators, namely
oscillator generators. For this reason we call the operators
(15) {\it preoscillator (creation and annihilation) generators} of 
$U[osp(2n+1/2m)]$ (resp. of $osp(2n+1/2m)$).
The preoscillator generators give an alternative to the Chevalley
description of $U[osp(2n+1/2m)]$. 

Observe that in this setting  the 
para-Bose (resp. the Bose) operators are odd, whereas the para-Fermi
(and the Fermi) operators are even generators. 

\smallskip
Coming back to the defining relations (16)
of the preoscillator generators
we note that they define a linear map
$$
L(n/m)\otimes L(n/m)\otimes L(n/m) \rightarrow L(n/m),\e(22)
$$
which identifies $osp(2n+1/2m)$ also as a Lie-supertriple system,
an approach which was recently developed in [18].

\s
Our purpose is to quantize $U[osp(2n+1/2m)]$ via the preoscillator
creation and annihilation operators. This is however difficult to be done
directly via the relations (16). Therefore in the next section we 
select a subset of relations from (16), which describe
completely $U[osp(2n+1/2m)]$, and which are convenient for quantization.

\vskip 6mm

\n{\bf 4. Description of $U[osp(2n+1/2m)]$ via Green generators} [13]
\bigskip

{\bf Proposition 3.} {\it  $U[osp(2n+1/2m)]$ is an associative unital
superalgebra with  generators 
$$
a_1^\pm,a_2^\pm,\ldots,a_{m-1}^\pm,a_{m}^\pm,a_{m+1}^\pm,
\ldots,a_{m+n}^\pm \equiv a_N^\pm,\e(23)
$$
referred as to Green generators,
relations $(\xi,\;\eta =\pm \; or \; \pm 1)$
$$
\eqalign{
& \[\[a_i^\eta,a_j^{-\eta}\],a_k^\eta\]=
2\eta^{\l k \r} \delta_{jk}a_i^\eta, \quad |i-j|\leq 1,
\q \eta =\pm , \cr
& [[a_{N-1}^\eta,a_{N}^{\eta}],a_{N}^\eta]=0, \q \eta =\pm , \cr
}\e(24)
$$
\n
and  ${\bf Z_2}$-grading
$$
{\rm deg}(a_i^\pm)= \langle i \rangle .\e(25)
$$
The Green generators (23) are the preoscillator generators of
$U[osp(2n+1/2m)]$.
}
\bigskip

In order to indicate how the proof can be done we recall the
Chevalley definition of $U[osp(2n+1/2m)]$ and write down explicit
relations between  the Green  and the Chevalley generators.  Let
$(\a_{ij}), \; i,j=1,\ldots,N$ be an $N\times N$ symmetric Cartan
matrix chosen as:
$$
(a_{ij})= (-1)^{\l j\r}
\delta_{i+1,j}+(-1)^{\l i \r}
\delta_{i,j+1} -[(-1)^{\l j+1
\r} + (-1)^{\l j \r}]\delta_{ij}+
\delta_{i,m+n}\delta_{j,m+n}. \eqno(26)
$$
For instance the Cartan matrix of $B(4/4) \equiv
osp(9/8)$ is $ 8 \times 8$ dimensional matrix:

\bigskip
$$
(a_{ij})=\pmatrix{
    2&-1& 0& 0& 0& 0& 0& 0\cr
   -1& 2&-1& 0& 0& 0& 0& 0\cr
    0&-1& 2&-1& 0& 0& 0& 0\cr
    0& 0&-1& 0& 1& 0& 0& 0\cr
    0& 0& 0& 1&-2& 1& 0& 0\cr
    0& 0& 0& 0& 1&-2& 1& 0\cr
    0& 0& 0& 0& 0& 1&-2& 1\cr
    0& 0& 0& 0& 0& 0& 1&-1\cr
}. \eqno(27)
$$

\bigskip

Then $U[osp(2n+1/2m)]$ is defined as an associative superalgebra
with 1 in terms of a number of generators subject to 
a number of relations.
The generators are the Chevalley generators $h_i,\;e_i,\;f_i,$
$i=1,\ldots,N$; the relations are the Cartan-Kac relations
$$
[h_i,h_j]=0,\;
[h_i,e_j]=a_{ij}e_j,\;
[h_i,f_j]=-a_{ij}f_j,\;
\[e_i,f_j\]=\delta_{ij}h_i, \eqno(28)
$$
the $e-$Serre relations
$$
\eqalign{
& [e_i,e_j]=0,\; {\rm for} \; \vert i-j \vert>1;\;\;
\[e_i,[e_i,e_{i\pm 1}]\]=0,\; i\neq N;\cr
& \{[e_{m-1},e_m],[e_m,e_{m+1}]\}=0 ;\;
[e_N,[e_N,[e_N,e_{N-1}]]]=0; \cr
}\e(29a)
$$
and the $f-$Serre relations
$$
\eqalign{
& [f_i,f_j]=0,\; {\rm for} \; \vert i-j \vert>1;\;\;
\[f_i,[f_i,f_{i\pm 1}]\]=0,\; i\neq N;\cr
& \{[f_{m-1},f_m],[f_m,f_{m+1}]\}=0 ;\;
[f_N,[f_N,[f_N,f_{N-1}]]]=0. \cr
}\e(29b)
$$
The grading on $U[osp(2n+1/2m)]$ is induced from:
$
{\rm deg}(e_m)={\rm deg}(f_m)=\bar 1,
\; {\rm deg}(e_i)={\rm deg}(f_i)=\bar 0
$ 
for $i\neq m.$

\smallskip
The expressions of the Green generators in terms of the Chevalley
generators read ($i=1,\ldots,N-1$):
$$
\eqalignno{
& a_i^- =(-1)^{(m-i)\l i \r}\sqrt{2}[e_i,[e_{i+1},
  [\ldots,[e_{N-2},[e_{N-1},e_N]]
  \ldots ]]], \q a_N^-=\sqrt{2}e_N, & \cr
&&\cr
& a_i^+=-\sqrt{2}[f_i,[f_{i+1},[\ldots,[f_{N-2},[f_{N-1},f_{N}]]
  \ldots ]]], \q a_N^+=-\sqrt{2}f_N. & (30)\cr  
}
$$

\n Then one proves that $a_i^\pm$ generate $U[osp(2n+1/2m)]$ 
($i=1,\ldots,N-1$),
$$
\eqalignno{
& h_i= {1\over 2}\[a_{i+1}^-,a_{i+1}^+\]-{1\over 2}\[a_i^-,a_i^+\],\q
  h_N=-{1\over 2}\[a_N^-,a_N^+\], & \cr
& e_i={1\over 2}\[a_i^-,a_{i+1}^+\], \q 
   e_N={1\over{\sqrt 2}}a_N^- ,& \cr
&  f_i={1\over 2}\[a_i^+,a_{i+1}^-\],\q
   f_N=-{1\over{\sqrt 2}}a_N^+ ,& (31)\cr
}
$$
and that the Cartan-Kac and the Serre relations follow from 
(24) and (31).

\vskip 6mm

\n{\bf 5. Description of $U_q[osp(2n+1/2m)]$ via deformed Green generators}

\bigskip
The $q$-deformed superalgebra $U_q[osp(2n+1/2m)]$, a Hopf superalgebra,
is by now a classical concept. See, for instance, [1-4] 
where all Hopf algebra operations are explicetly given. Here,
following [4], we write only the algebra operations.

{\bf Proposition 4.} {\it $U_q[osp(2n+1/2m)]$ is an associative unital 
algebra with Chevalley generators \hfill\break
$ e_i,\; f_i,\; k_i=q^{h_i},\; {\bar k}_i \equiv k_i^{-1}=q^{-h_i},\; 
i=1,\ldots,N, $  which satisfy the Cartan-Kac relations 
$$
\eqalign{
&  k_ik_i^{-1}=k_i^{-1}k_i=1, \q  k_ik_j=k_jk_i,\cr
&  k_ie_j=q^{\alpha_{ij}}e_jk_i,\q k_if_j=q^{-\alpha_{ij}}f_jk_i,\cr
&  \[e_i,f_j\]=\delta_{ij}{{k_i-{\bar k}_i}\over{q-{\bar q}}}, 
}\eqno(32)
$$
the $e-$Serre relations 
$$
\eqalign{
& (e1)\q \[e_i,e_j\]=0, \quad \vert i-j \vert \neq 1, \cr
& (e2)\q [e_i,[e_i,e_{i \pm 1}]_{\bar q}]_q \equiv
[e_i,[e_i,e_{i \pm 1}]_q]_{\bar q}=0, \quad
 i\neq m,\;\; i\neq N \cr 
& (e3)\q  \{[e_{m},e_{m-1}]_{q},[e_m,e_{m+1}]_{\bar q}\}=0,\cr
& (e4)\q  [e_N,[e_N,[e_N,e_{N-1}]_{\bar q}]]_q \equiv
[e_N,[e_N,[e_N,e_{N-1}]_q]]_{\bar q}=0, \cr 
}\eqno(33)
$$
and the $f-$Serre relations
$$
\eqalign{
& (f1)\q \[f_i,f_j\]=0, \quad \vert i-j \vert \neq 1, \cr
& (f2)\q [f_i,[f_i,f_{i \pm 1}]_{\bar q}]_q \equiv
[f_i,[f_i,f_{i \pm 1}]_q]_{\bar q}=0, \quad
 i\neq m,\;\; i\neq N \cr 
& (f3)\q \{[f_{m},f_{m-1}]_{q},[f_m,f_{m+1}]_{\bar q}\}=0,\cr
& (f4)\q  [f_N,[f_N,[f_N,f_{N-1}]_{\bar q}]]_q \equiv
[f_N,[f_N,[f_N,f_{N-1}]_q]]_{\bar q}=0. \cr 
}\eqno(34)
$$
}
The $(e3)$ and $(f3)$ Serre relations are the additional Serre
relations [19-21], which were initially omitted.

\s
We are now ready to state our main result, given also in the 
Introduction.

\s
{\bf Theorem.} {\it $U_q[osp(2n+1/2m)]$ is an associative 
superalgebra with 1, generators 
$a_i^\pm,\; L_i,\;{\bar L_i}\equiv L_i^{-1}$, $i=1,2,\ldots,m+n=N$,
relations ($\xi,\eta =\pm\; or\; \pm 1,\;\;{\bar q}\equiv q^{-1}$)
$$
\eqalign{
&  L_iL_i^{-1}=L_i^{-1}L_i=1, \q  L_iL_j=L_jL_i,\cr
& L_ia_j^\pm=q^{\pm \delta_{ij}(-1)^{\l i \r} }a_j^\pm,\cr
& \[a_i^-,a_i^+\]=-2{L_i-{\bar L}_i\over q-{\bar q}}, \cr 
& \[\[a_i^{\eta},a_{i+ \xi}^{-\eta}\],  
a_j^{\eta}\]_{q^{-\xi (-1)^{\l i \r}\delta_{ij}}}
=2(\eta)^{\l j \r}\delta _{j,i + \xi}L_j^{-\xi \eta}a_i^{\eta}, \cr  
& [[a_{N-1}^\xi,a_N^\xi],a_{N}^\xi]_{\bar q}=0 . \cr
}\eqno(35)
$$
and ${\bf Z}_2$-grading 
$ {\rm deg}(L_i)={\bar 0},\;{\rm deg}(a_i^\pm)= \langle i \rangle $. 
}

\bigskip
The expressions  of $a_i^\pm$ and $L_i$ via the Chevalley generators
read ($i=1,\ldots,N-1$):
$$
\eqalign{
&  L_i=k_ik_{i+1}\ldots k_N \;({\rm including} \;i=N),  \cr
&  a_i^-=(-1)^{(m-i)\l i \r}\sqrt{2}[e_i,[e_{i+1},
  [\ldots,[e_{N-2},[e_{N-1},e_N]_{q_{N-1}}]_{q_{N-2}}
  \ldots ]_{q_{i+2}}]_{q_{i+1}}]_{q_{i}},\q 
   a_N^-=\sqrt{2}e_N,   \cr
&  a_i^+=(-1)^{N-i+1}\sqrt{2}
   [[[\ldots [f_N,f_{N-1}]_{{\bar q}_{N-1}},f_{N-2}]_{{\bar q}_{N-2}} 
   \ldots]_{{\bar q}_{i+2}},f_{i+1}]_{{\bar q}_{i+1}},f_{i}]_{{\bar q}_{i}}
   , \q a_N^+=-\sqrt{2}f_N,  \cr  
}\eqno(36)
$$
\n where 
$$
q_i={\bar q},\;i=1,\ldots,m-1;\q q_i=q,\;i=m,\ldots,N. 
$$

\s
The next result is essential for the proof of the Theorem.

\vfill\eject
{\bf Proposition 5.} {\it The following relations hold:}
$$
\eqalignno{
1.\q & \[e_i,a_j^+\]=-\delta_{ij}(-1)^{\l i+1 \r}k_ia_{i+1}^+,
       \q  i\ne N, & (37)\cr
&&\cr
2.\q & \[a_j,f_i\]=\delta_{ij}a_{i+1}^-{\bar k}_i,\q i\ne N, & (38)\cr
&&\cr
3.\q & \[e_i,a_j^-\]=0, \;\;if\; i<j-1 \; or \; i>j,\q i\ne N,& (39a)\cr
     & \[e_i,a_{i+1}^-\]_{q_i}=(-1)^{\l i+1 \r}a_i^-, 
       \q i\ne N,   &  (39b) \cr
     & \[e_i,a_{i}^-\]_{{\bar q}_{i-1}}=0,\q i\ne N, & (39c) \cr
&&\cr
4.\q & \[a_j^+,f_i\]=0, \;\;if\; i<j-1 \; or \; i>j,\q i\ne N,& (40a)\cr
     & \[a_{i+1}^+,f_i\]_{{\bar q}_i}=-a_i^+,\q i\ne N. & (40b) \cr
     &  \[a_i^+,f_i\]_{q_{i-1}}=0, \q i\ne N. & (40c) \cr
}
$$

Also here one proves that $a_i^\pm$ and $L_i^{\pm 1}$ 
generate $U_q[osp(2n+1/2m)]$. More precisely ($i=1,\ldots,N-1$),
$$
\eqalignno{
&  k_i= L_i{\bar L}_{i+1},\q L_N=k_N, & \cr
&  e_i={1\over 2}{\bar L}_{i+1}\[a_i^-,a_{i+1}^+\], \q 
   e_N={1\over{\sqrt 2}}a_N^- ,& \cr
&  f_i={1\over 2}\[a_i^+,a_{i+1}^-\]L_{i+1},\q
   f_N=-{1\over{\sqrt 2}}a_N^+ .& (41)\cr
}
$$

It is a long computation to show, using only the relations (35),
that the operators (41) satisfy the Cartan-Kac and the Serre relations. 
The proof is based on repeated use of nontrivial identities. 
Here is one of them.

\s
{\bf Proposition 6.} {\it
If $B$ or $C$ is an even element, then for any values of the parameters
$x,y,z,t,r,s$ subject to the relations
$$
x=zs,\;\; y=zr,\;\; t=zsr, \eqno(42)
$$
the following identity holds:
$$
\[A,[B,C]_x\]_y=\[\[A,B\]_z,C\]_t +
(-1)^{deg(A)deg(B)}z\[B,\[A,C\]_r\]_s.\eqno(43)
$$
}

In particular it is nontrivial to prove that $e_m^2=0$, which is one
of the Serre relations, or to show that the additional Serre
relations $(e3)$ and $(f3)$ hold.

\vskip 6mm

{\bf 5. Concluding remarks}

\bigskip
The root system of the orthosymplectic Lie superalgebra
$osp(2n+1/2m)$ reads:
$$
\Delta=\{\xi \varepsilon_i +\eta \varepsilon_j; \xi \varepsilon_i;
2\xi \varepsilon_k \},\q
i\neq j=1,\ldots,m+n\equiv N ;\;k=1,\ldots,m;
\; \xi, \eta =\pm\}.
\eqno(44)
$$
The roots $\varepsilon_1,\ldots,\varepsilon_N$ are orthogonal
with respect to the Killing form on $osp(2n+1/2m)$.
The Green generators are the root vectors, corresponding to the
orthogonal roots. More precisely, the correspondence reads:
$$
a_i^\pm \leftrightarrow \mp \varepsilon_i, \q
i=1,\ldots,N. \eqno(45)
$$
Therefore what we have done here is 

\s
(1) to describe  $U[osp(2n+1/2m)]$
in terms of a ``minimal'' set of relations among the positive and
the negative root vectors, corresponding to the orthogonal roots.

\s
(2) to describe  $U_q[osp(2n+1/2m)]$ 
entirely in terms of deformed ``orthogonal'' root vectors,
namely deformed Green generators.

\s
This is a good opportunity to mention that the canonical quantum statistics
and its generalization, the parastatistics, is based on the
representation theory of orthosymplectic Lie superalgebras.
For instance the Bose operators $B_i^\pm,\;i=1,\ldots,n$ are generators
of $osp(1/2n)$ in a particular representation. Similar statement
holds for $n$ pairs of Fermi creation and annihilation operators:
they are generators of the Lie algebra $so(2n+1)$ in a particular, the
Fock representation. Both $osp(1/2n)$ and $so(2n+1)$
are among the superalgebras from the class {\it B}
in the classification of Kac of the basic Lie superalgebras [14].
Therefore the canonical quantum statistics and its generalization,
the parastatistics, could be called {\it B-statistics}.

\s
One can associate a concept of creation and annihilation operators
with every simple Lie algebra [22-24] and presumably also
with every basic Lie superalgebra. The creation and 
the annihilation operators of the Lie superalgebra sl(1/n) were
given in [24]. Therefore, parallel to the {\it B-statistics}, 
i.e., the parastatistics, there exists {\it A-statistics},  
{\it C-statistics} and {\it D-statistics}. The corresponding
deformations, certainly, also exist. 
In fact the {\it A-statistics} belongs to the
class of the exclusion statistics, recently introduced by
Haldane [25] in solid state physics.
\bigskip

\vskip 20pt
\noindent
{\it Acknowledgments.}
I am thankful to Prof. Doebner for the invitation to report the
present investigation at the XXI Colloquium on Group Theoretical
Methods in Physics. It is a pleasure
to thank Prof. C. Reina for making it possible for me to visit the
Mathematical Physics Sector in Sissa,  where most of the
results were obtained and Prof. Randjbar-Daemi for the kind
hospitality at the High Energy Section of ICTP.
The work was supported also by the Grant
$\Phi-416$ of the Bulgarian Foundation for Scientific Research.

\vskip 24pt
\noindent
{\bf References}

\bigskip
\settabs \+  [11] & I. Patera, T. D. Palev, Theoretical 
   interpretation of the experiments on the elastic \cr

\+ [1] & Chaichian M. and Kulish P. 1990 {\it Phys. Lett.} {\bf 234B} 72\cr

\+ [2] & Bracken A.J., Gould M.D. and Zhang R.B. 1990 {\it Mod. Phys. Lett. A}
         {\bf 5} 331 \cr

\+ [3] & Floreanini R., Spiridonov V.P. and Vinet L. 1991 
         {\it Comm. Math. Phys.} {\bf 137} 149 \cr

\+ [4] & Khoroshkin S.M. and Tolstoy V.N. 1991
         {\it Comm. Math. Phys.} {\bf 141} 599 \cr

\+ [5] & Green H.S. 1953 {\it Phys. Rev.} {\bf 90} 270 \cr

\+ [6] & Floreanini R. and Vinet L. 1990 {\it J. Phys. A} {\bf 23}
         L1019  \cr

\+ [7] & Celeghini E., Palev T.D. and Tarlini M. 1991 
        {\it Mod. Phys. Lett. B } {\bf 5} 187 \cr
   
\+ [8] & Palev T.D. 1993 {\it J. Phys. A} {\bf 26} L111 \cr

\+ [9] & Hadjiivanov L.K. 1993 {\it J. Math. Phys.} {\bf 34} 5476 \cr

\+ [10] & Palev T.D. and Van der Jeugt J.  1995 {\it J. Phys. A} 
          {\bf 28} 2605 \cr

\+ [11] & Palev T.D. 1994 {\it Lett. Math. Phys.} {\bf 31} 151 \cr

\+ [12] & Palev T.D. 1993 {\it J. Math. Phys.} {\bf 34} 4872 \cr

\+ [13] & Palev T.D. 1996 {\it J. Phys. A} {\bf 29} L171 \cr

\+ [14] & Kac V.G. 1978 {\it Lect.Notes Math.} {\bf 676} 597 \cr

\+ [15] & Kamefuchi S. and Takahashi Y. 1960 {\it Nucl. Phys.}
         {\bf 36} 177\cr

\+ [16] & Omote M., Ohnuki Y. and Kamefuchi S.
         1976 {\it Prog. Theor. Phys.} {\bf 56} 1948\cr

\+ [17] & Ganchev A. and Palev T.D. 1980 {\it J. Math. Phys.} {\bf 21}
         797 \cr

\+ [18] & Okubo S. 1994 {\it J. Math. Phys.} {\bf 35} 2785 \cr

\+ [19] & Khoroshkin S.M. and Tolstoy V.N. 1991 {\it Comm. Math. Phys.}
        {\bf 141} 599 \cr

\+ [20] & Floreanini R., Leites D. A. and Vinet L. 1991		
         {\it Lett. Math. Phys.} {\bf 23 } 127\cr

\+ [21] & Scheunert M. 1992 {\it Lett. Math. Phys.} {\bf 24 } 173 \cr

\+ [22] & Palev T.D. 1976 Thesis (Institute of Nulear Research 
          and Nuclear Energy, Sofia)\cr

\+ [23] & Palev T.D. 1979 {\it Czech. J. Phys.} {\bf B29} 91\cr

\+ [24] & Palev T.D. 1980 {\it J. Math. Phys.} {\bf 21} 1293 \cr

\+ [25] & Haldane F.D. 1991 {\it Phys. Rev. Lett.} {\bf 67} 937\cr

\end